\documentclass[12pt]{amsart}

\usepackage{amsthm}
\usepackage{amsmath,amssymb}
\usepackage{geometry}

\begin{document}

\newtheorem{theorem}{Theorem}
\newtheorem{proposition}[theorem]{Proposition}
\newtheorem{corollary}[theorem]{Corollary}
\newtheorem{lemma}[theorem]{Lemma}
\newtheorem{definition}[theorem]{Definition}
\newtheorem{remark}[theorem]{Remark}
\newtheorem{conjecture}[theorem]{Conjecture}
\newtheorem{example}[theorem]{Example}

\newcommand\Z{{\mathbb Z}}
\newcommand\K{{\mathbb K}}
\newcommand\R{{\mathbb R}}
\newcommand\A{{\mathcal A}}
\newcommand\D{{\mathcal D}}
\newcommand\C{{\mathcal C}}
\newcommand\AV{{\mathcal {AV}}}
\newcommand\F{{\mathcal F}}
\newcommand\V{{\mathcal V}}
\newcommand\LL{{\mathcal L}}
\newcommand\hL{\widehat{\mathcal L}}
\newcommand\dd[1]{\frac{\partial}{\partial #1}}
\newcommand\Der{\textup{Der\,}}
\newcommand\Hom{\textup{Hom}}
\newcommand\Span{\textup{Span}}
\newcommand\End{\textup{End}}
\newcommand\ad{\textup{ad}}
\newcommand\tr{\textup{tr}}
\newcommand\id{\textup{id}}
\newcommand\gl{\mathfrak{gl}}
\newcommand\g{\mathfrak{g}}
\newcommand\m{\mathfrak{m}}
\newcommand\sll{\mathfrak{sl}}
\newcommand\NN{{\mathcal N}}

\title
[$AV$-modules on affine algebraic groups]
{$AV$-modules on affine algebraic groups}
\author{Yuly Billig}
\address{School of Mathematics and Statistics, Carleton University, Ottawa, Canada}
\email{billig@math.carleton.ca}
\author{Colin Ingalls}
\email{ColinIngalls@carleton.ca}
\let\thefootnote\relax

\maketitle

\begin{abstract}
We study $AV$-modules on an affine algebraic group $G$, that is, modules that admit compatible actions of the algebra of regular functions on $G$ and the Lie algebra of vector fields. A central role in this theory is played by the Lie algebra of jets of vector fields. Our main result is the Decomposition Theorem for the Lie algebra of jets of vector fields. We apply the Decomposition Theorem to describe differentiable $AV$-modules on $G$.
\end{abstract}

\section{Introduction}

The theory of $D$-modules on algebraic groups has been instrumental in the proof of the Kazhdan-Lusztig conjecture by Brylinski-Kashiwara and Beilinson-Bernstein \cite{BK, BB, H}. In this paper, we discuss $AV$-modules on algebraic groups, which are natural generalizations of $D$-modules. Geometric objects like tensor fields do not form a $D$-module under the Lie derivative action, but belong to the category of $AV$-modules.
The theory of $AV$-modules was successfully used in establishing classifications of simple Harish-Chandra modules for Lie algebras of vector fields on a torus \cite{BF0}, and on the affine space \cite{GS, XL}.

An $AV$-module on a connected affine algebraic group $G$ is a module that admits an action of the algebra $A$ of regular functions on $G$, and an action of the Lie algebra $V$ of vector fields, with the Leibniz rule as a compatibility condition (see Definition \ref{AV} in the next Section). The pair of algebras $(A, V)$ forms a Rinehart pair 
\cite{Ri}, and the category of $AV$-modules is governed by the weak enveloping algebra $U_\text{weak}(A, V)$, which is the smash product of the Hopf algebra $U(V)$ with its module $A$:  $U_\text{weak}(A, V)
= A \# U(V)$. This algebra was first introduced by Jacobson in \cite{Jac}. This associative algebra has a Lie subalgebra $A\#V$, the smash product of the Lie algebra $V$ with a commutative algebra $A$, on which $V$ acts by derivations. The Lie algebra $A\#V$ has a geometric interpretation as the Lie algebra of algebraic jets of vector fields, and its completion, $A\widehat\# V$, is the Lie algebra of $\infty$-jets of vector fields. 

The Lie algebra $A\widehat\# V$ has an interesting structure, and is important for understanding the theory of $AV$-modules. The central result of the paper is the Decomposition Theorem (Theorem \ref{decomposition} in Section 4), which states that on an affine algebraic group $G$, there exists a semidirect product decomposition 
$$A\widehat\# V \cong V \ltimes (A \widehat\otimes \LL_+),$$
where $\LL_+$ is the Lie algebra of vector fields on $G$, vanishing at the identity element $e$. 
We point out that this is a non-trivial isomorphism -- in the semidirect product, neither $V$ is mapped to the Lie algebra of vector fields $1 \# V \subset A\widehat\# V $, nor $\LL_+$ is identified with a subalgebra in $1 \# V$. 
Previously, such a Decomposition Theorem was known for affine space \cite{XL, BIN}, the torus, and only locally for an affine variety \cite{BI}.

In order to prove the Decomposition Theorem, we generalize the classical Taylor formula 
$$f(x) = \sum_{k=0}^\infty \frac{1}{k!} \,\frac{d^k f}{d x^k} (a) \,(x-a)^k$$
to the setting of functions on algebraic groups (Theorem \ref{Taylor} in Section 3). We replace the family
$\left\{\left(\frac{d}{dx}\right)^k\right\}$ with a basis of the universal enveloping algebra $U(\g^L)$ of the Lie algebra $\g^L$ of left-invariant vector fields, and view $\left\{ (x-a)^k \right\}$ as a basis for the algebra of $G$-invariants in $A \otimes A$.  

The most important class of $AV$-modules is the category of differentiable $AV$-modules. An $AV$-module $M$ is called $N$-differentiable if it is annihilated by the elements 
$$\sum\limits_{k=0}^\infty (-1)^k {N \choose k} f^{N-k} \# f^k \tau \ \ \text{for all \, } f\in A, \tau\in V.$$
$D$-modules are precisely 1-differentiable $AV$-modules, and all tensor modules are 2-differentiable.
Using the Decomposition Theorem, we prove that $N$-differentiable $AV$-modules on $G$ are modules for the associative algebra $A \# U(\g \ltimes (A \otimes (\LL_+ / \m_e^{N-1} \LL_+)))$,
where $\g$ is the Lie algebra of $G$, and $\m_e$ is the ideal of functions on $G$ vanishing at $e$.

Our description is more detailed for 2-differentiable $AV$-modules on $G$. Such a module is given by a $D$-module structure on $M$, together with an action of $\gl_n$, $n = \dim G$, which commutes with the action of $D$. In Theorem  \ref{twodiff}, we give an explicit formula for this correspondence.

In the final section of this paper, we point out that the Lax equation emerges in the context of 2-differentiable $AV$-modules, connecting them to the theory of integrable dynamical systems.

The paper is organized as follows. In Section 2, we outline the basic definitions related to $AV$-modules, Rinehart pairs, and jets of vector fields. We generalize the Taylor formula to the context of affine algebraic groups in Section 3. In Section 4, we establish the Decomposition Theorem, the main result of the paper. In Section 5, we present examples of $AV$-modules, focusing on 2-differentiable modules.

\subsection*{Acknowledgements} Both authors gratefully acknowledge support from the Natural Sciences and Engineering Research Council of Canada through Discovery Grants. We benefited from helpful discussions with ChatGPT.

\section{$AV$-modules and jets of vector fields}

Let $A$ be a commutative associative algebra over a field $K$ of characteristic zero.
Assume that the Lie algebra $V$ acts on $A$ by derivations.

\begin{definition}
\label{AV}
An $AV$-module $M$ is a module for $A$ and also a module for the Lie algebra $V$, with two actions compatible via the Leibniz rule:
\begin{equation}
\label{Leibniz}
\tau (f m) = \tau(f) m + f ( \tau m) \quad \text{for \ } \tau\in V, \, f \in A, \, m \in M.
\end{equation}
\end{definition} 

Just as the universal enveloping algebra $U(L)$ controls the category of modules for Lie algebra $L$, there is an associative algebra that controls the category of $AV$-modules. The weak enveloping algebra $U_{\text{weak}}(A, V) = A \# U(V)$ is the smash product of the Hopf algebra $U(V)$ with its module $A$ \cite{Jac}. As a vector space, it is the tensor product $A \otimes_K U(V)$ with multiplication 
$$(f \otimes u) (h \otimes v) = \sum\limits_i f u_i^1 (h) \otimes u_i^2 v,$$
where $u \mapsto  \sum\limits_i  u_i^1 \otimes  u_i^2$ is the coproduct written in Sweedler notations.
An $AV$-module is the same as a module over the associative algebra $A \# U(V)$.

It turns out that the subspace $A\#V$ is a Lie subalgebra of $A \# U(V)$.

\begin{definition}
Let $A$ be a commutative associative algebra on which the Lie algebra $V$ acts by derivations. The smash product $A \# V$ is a Lie algebra whose underlying space is $A \otimes V$ and the Lie bracket is given by
\begin{equation}
\label{smash}
[f \# \tau, h \# \eta] = f \tau(h) \# \eta - h \eta(f) \# \tau + fh \# [\tau, \eta] \ \ \text{for\,}
f, h \in A, \, \tau,\eta \in V.
\end{equation}
\end{definition}

We are particularly interested in the case where $A$ is the algebra of regular functions on a smooth irreducible affine algebraic variety $X$ and $V = \Der(A)$.  In this setting, $V$ is the Lie algebra of vector fields on $X$.
Such a pair $(A, V)$ is an example of a Rinehart pair \cite{Ri}:

\begin{definition}
Let $A$ be a commutative associative algebra, and let $V$ be a Lie algebra. The pair $(A, V)$ is a Rinehart pair if $V$ acts on $A$ by derivations, and $V$ is an $A$-module, satisfying:
\begin{align*}
[\tau, f\eta] &= \tau(f) \eta + f [\tau,\eta], \\
(f\eta)(h) &= f (\eta(h)) \quad \text{for \ }  f, h \in A,\, \tau, \eta \in V.
\end{align*}
\end{definition}

There are two natural categories of modules for a Rinehart pair: $AV$-modules, also known as weak modules,  and $D$-modules, also known as strong modules, satisfying (\ref{Leibniz}) as well as the following axiom:
\begin{equation}
\label{Dmod}
(f \tau) m = f (\tau m), \ \ \text{for\,} f \in A, \, \tau \in V, \, m \in M.
\end{equation}

The algebra $A$ itself is a $D$-module for a Rinehart pair $(A, V)$, while $V$ is an $AV$-module, but not a $D$-module in general. 

The category of $D$-modules is controlled by the strong enveloping algebra $U_{\text{strong}}(A, V)$
\cite{Ri}, which  is the quotient
$$A \# U(V) / \left< fh \# \tau - f\#h\tau \ \big| \, f, h \in A, \tau\in V \right>.$$

When $A$ is the algebra of regular functions on a smooth irreducible affine variety $X$, and $V = \Der (A)$, the strong enveloping algebra $U_{\text{strong}}(A, V)$ is the algebra $D$ of differential operators on $X$ 
\cite{Ba}.
As we shall see below, in this context, the Lie algebra $A \# V$ may be interpreted as the Lie algebra of algebraic jets of vector fields on $X$.

If $V$ acts on $A$ by derivations, then the pair $(A, A\# V)$ is a Rinehart pair. Moreover,
$$U_{\text{strong}}(A, A\#V) = U_{\text{weak}}(A, V).$$

Let $\Delta \subset A \otimes A$ be the kernel of the multiplication map $\mu:\,A \otimes A \rightarrow A$. The quotient $A \otimes A/ \Delta^{N+1}$ is the algebra of $N$-jets of functions on $X$, and the completion
$$A \widehat\otimes A = \varprojlim\limits_N A \otimes A / \Delta^N$$
is the algebra of $\infty$-jets of functions (see Chapter 10 in \cite{AM} for the basics on the completions of algebras).

The space $\Delta^N \otimes_A V$ is an ideal in the Lie algebra $A \# V$, and the quotient  \break
$A \# V/ \Delta^{N+1} \otimes_A V$ is the Lie algebra of $N$-jets of vector fields on $X$. The completion
$$ A \widehat\# V =  \varprojlim\limits_N (A \# V) / (\Delta^N \otimes_A V)$$
is the Lie algebra of jets of vector fields on $X$ \cite{BI}. Since $\mathop\cap\limits_{N=1}^\infty \Delta^N = (0)$, the Lie algebra $A\#V$ injects into the Lie algebra of $\infty$-jets of vector fields. Thus, $A\#V$ may be interpreted as the Lie algebra of algebraic jets of vector fields.

We now turn to the principal setting of this paper, the case when the underlying variety is a connected affine algebraic group, $X = G$.

There are three natural actions of $G$ on itself: left translation: $L_g (h) = gh$, right translation: $R_g (h) = hg^{-1}$, and the action by conjugation: $C_g (h) = g h g^{-1}$, with $g, h \in G$. 
Note that the left and right actions of $G$ commute.
These induce three actions of $G$ on the algebra $A$ of regular functions on $G$:  $L^*_g (f) = f \circ L_{g^{-1}}$, $R^*_g (f) = f \circ R_{g^{-1}}$, $C^*_g (f) = f \circ C_{g^{-1}}$  for $f \in A$. 
Likewise, there are three actions of $G$ on $V = \Der A$: $L^\circ_g \tau = L^*_g \circ \tau \circ L^*_{g^{-1}}$, etc. A vector field $\tau \in V$ is called left-invariant if $L^\circ_g \tau = \tau$ for all $g \in G$. 

Let $\g$ be the tangent space to $G$ at the identity $e$, which may be identified with the space of derivations $A \rightarrow K$ at the identity: $x(fh) = x(f) h(e) + f(e) x(h)$, for $x \in \g$, $f, h \in A$. We construct a left-invariant vector field $x^L$ which restricts to $x$ at the identity: $x^L (f) (g) = x(f \circ L_g)$. Likewise, we construct a right-invariant vector field $x^R$, which restricts to $-x$ at the identity: $x^R (f) (g) = - x(f \circ R_{g^{-1}})$. The negative sign is needed in order for the two constructions to induce the same Lie bracket on $\g$. We get two isomorphic subalgebras $\g^L$ and $\g^R$ of left- and right-invariant vector fields in $V$. Alternatively, $x^L$ and $x^R$ may be expressed using the coproduct $A \rightarrow A \otimes A$, $f \mapsto \sum_i f^1_i \otimes f^2_i$. Then $x^L(f) = \sum_i f^1_i  x(f^2_i)$, $x^R(f) = -\sum_i x(f^1_i)  f^2_i$.

Note that the elements of $\g^L$ are the differentials of the {\emph right} action of $G$ on $A$, and $\g^R$ are the differentials of the {\emph left} action.

There is also a homomorphism $\g \rightarrow V$ induced by the conjugation action, $x \mapsto x^C = x^L + x^R$. The kernel of this homomorphism is $Z(\g)$. Since $x^C$ evaluates to zero at the identity, it preserves the maximal ideal $\m_e$.

The Lie algebra of vector fields on an affine algebraic group $G$ is a free $A$-module generated by either $\g^L$ or $\g^R$ \cite{BF}:
$$V = A \otimes \g^L = A \otimes \g^R.$$ 
In fact, these tensor products are smash products of $A$ with $\g^L$ and $\g^R$. For the Lie algebra of algebraic jets of vector fields, we have the double smash product construction:
$$A \# V = A \# (A \# \g^L).$$
In this setting, Lie bracket takes the following form:
\begin{align}
\label{doublesmash}
[p \otimes x^L, q \otimes y^L] =& (\mu(p) \otimes 1) (x^L \otimes 1)(q) \otimes y^L
+ p (1 \otimes x^L)(q)\otimes y^L \\
&- (\mu(q) \otimes 1) (y^L \otimes 1)(p) \otimes x^L 
- q (1 \otimes y^L)(p)\otimes x^L \nonumber\\
&+ pq \otimes [x,y]^L, \  \ \text{for } p,q \in A\otimes A, x,y \in \g. \nonumber
\end{align}

\section{Taylor formula for functions on algebraic groups}

In this section, we will establish the Taylor formula for functions on algebraic groups.

The algebra $A \otimes A$ of functions on $G \times G$ is a $G$-module under the diagonal left action of $G$. Denote by $(A \otimes A)^L$ the subalgebra of $G$-invariants under this action. Since the differential of the left action of $G$ on $A$ yields the action of the Lie algebra of right-invariant vector fields $\g^R$, we conclude that $(A \otimes A)^L$ is the kernel of the action of $\g^R$ on the tensor product $A \otimes A$.

\begin{lemma}
$(A \otimes A)^L \cong A$ as algebras. 
\end{lemma}
\begin{proof}
The isomorphism maps $\varphi: (A \otimes A)^L \rightarrow A$ and $\psi: A \rightarrow  (A \otimes A)^L$ are given by
$$\varphi\left( \sum_i f_i \otimes h_i \right) (g) = \sum_i f_i (e) h_i (g) \text{\ \ and \ }
\psi(f) (g_1, g_2) = f(g_1^{-1} g_2).$$  
It is easy to see that the image of $\psi$ is consists of left-invariant functions, and that these maps are the inverses of each other:
\begin{align*}
&\varphi(\psi(f)) (g) = \psi(f) (e, g) = f(e^{-1} g) = f(g), \\
&\psi\left( \varphi\left( \sum_i f_i \otimes h_i \right) \right) (g_1, g_2)
= \varphi\left( \sum_i f_i \otimes h_i \right) (g_1^{-1} g_2) = \\
&= \sum_i f_i (e) h_i (g_1^{-1} g_2)
= \left( \sum_i f_i \otimes h_i \right) (g_1, g_2),
\end{align*}
where in the last equality follows from left-invariance.
\end{proof}

Since the left and right actions of $G$ commute, $(A \otimes A)^L$ is a $G$-module under the actions $R^* \otimes \id$, $\id \otimes R^*$, and  $R^* \otimes R^*$. It is straightforward to check the following claims:

\begin{lemma}
Under the isomorphism $(A \otimes A)^L \cong A$, the $R^* \otimes \id$ action corresponds to the left $G$-action $L^*$ on $A$, $\id \otimes R^*$ corresponds to the right $G$-action $R^*$, and $R^* \otimes R^*$ corresponds to the conjugation action $C^*$ of $G$ on $A$. 
\end{lemma} 
\begin{corollary}
\label{transfer}
Under the isomorphism $\varphi: (A \otimes A)^L \rightarrow A$,
 (i) $\varphi \circ (x^L \otimes 1) = x^R \circ \varphi$, 
(ii) $\varphi \circ (1 \otimes x^L) = x^L \circ \varphi$, 
(iii) $\varphi \circ (x^L \otimes 1 + 1 \otimes x^L) = x^C \circ \varphi$,
for $x \in \g$.  
\end{corollary}
\begin{proof}
This follows from the fact that left-invariant vector fields are the differentials of the $R^*$-action, and right-invariant vector fields are the differentials of the $L^*$-action of $G$.
\end{proof}

Consider the pairing $U(\g^L) \times A \rightarrow K$, given by applying the differential operator $u \in U(\g^L)$ to function $f \in A$ and evaluating the result at the identity element of $G$.

Let $U_0 \subset U_1 \subset U_2 \subset \ldots$ be  the standard augmentation filtration on $U(\g^L)$. Here $U_0 = K \cdot 1$, $U_1 = K \cdot 1 + \g^L$, etc.
Consider also the decreasing filtration on $A$ by powers of the maximal ideal $\m_e$:
$A \supset \m_e \supset \m_e^2 \supset \ldots$. Since $U_n (\g^L)$ acts on $A$ by differential operators of order less or equal to $n$, it trivially pairs to $\m_e^{n+1}$, and we get a pairing of finite-dimensional spaces:
$$U_n(\g^L) \times A/\m_e^{n+1} \rightarrow K.$$
By Cartier's Theorem (\cite{DG} II.6.1), this pairing is non-degenerate, and we can choose dual bases with respect to this pairing. By increasing the value of $n$, inductively, we construct dual bases $\{ u_i \}$ in $U(\g^L)$ and $\{ F_i \}$ in $A$. 

Recall the classical Taylor formula:
$$f(x) = \sum_{k=0}^\infty \frac{1}{k!} \,\frac{d^k f}{d x^k} (a) \,(x-a)^k.$$
We view both sides as functions on $\R \times \R$, where the first copy is parametrized by $a$ and the second copy by $x$. Note that the expression $x-a \in A \otimes A$ is translation invariant (under the diagonal action of the additive group $\R$ on $\R \times \R$). Thus, $(x-a)^k \in (A \otimes A)^L$. The Taylor formula yields the statement $1 \otimes A \subset A \widehat\otimes (A \otimes A)^L$. The following theorem generalizes this to an arbitrary affine algebraic group. 

\begin{theorem} [Taylor formula]
\label{Taylor}
Let $A$ be the algebra of regular functions on an affine algebraic group $G$ over a field $K$ of characteristic zero. Then
$$A \widehat\otimes A \cong A \widehat\otimes (A \otimes A)^L,$$
given by the Taylor formula
$$f(g_1, g_2) = \sum_i ( (1 \otimes u_i) f) (g_1, g_1) \,F_i (g_1^{-1} g_2).$$
The completion in the left hand side is
$$A \widehat\otimes A = \varprojlim_N A \otimes A/ \Delta^N,$$
and the completion in the right hand side is
$$A \widehat\otimes (A \otimes A)^L = \varprojlim_N A \otimes (A \otimes A)^L / A \otimes \psi(\m_e^N).$$
Applying the isomorphism $(A \otimes A)^L \cong A$, we obtain a non-trivial isomorphism
\begin{equation}
\label{transform}
A \widehat{\mathop\otimes\limits_\Delta} A \cong A \mathop{\widehat\otimes}\limits_{1 \otimes \m_e} A.
\end{equation}
The above maps are isomorphisms of $A$-modules with respect to the action of $A$ on the left tensor factors.
\end{theorem}
\begin{proof}
First, let us establish the formula
\begin{equation}
\label{T}
f(g_2) = \sum_i u_i (f) (g_1) F_i(g_1^{-1} g_2).
\end{equation}
Cartier's pairing yields a generalization of the Maclaurin formula:
$$f(g) = \sum_i u_i (f) (e) F_i (g)$$
in the completion $\widehat{A}$ of $A$ with respect to the $\m_e$-adic topology. Applying the Maclaurin
formula to $L^*_{g_1^{-1}} f$, we get
$$(L^*_{g_1^{-1}}f)(g) = \sum_i u_i (L^*_{g_1^{-1}} f) (e) F_i (g).$$
Since $u_i$ are left-invariant, this yields
$$f( g_1 g) = \sum_i u_i (f) (g_1) F_i(g).$$
Setting $g = g_1^{-1} g_2$, we obtain (\ref{T}).
Next, we apply (\ref{T}) to the second argument of $f(g_1, g_2)$:
$$f(g_1, g_2) = \sum_i ((1 \otimes u_i) (f)) (g_1, g_1) F_i(g_1^{-1} g_2).$$
Note that the change of variables $g = g_1^{-1} g_2$ transforms $\m_e$-adic topology into $\Delta$-adic topology, $\psi(\m_e) = \Delta \cap (A \otimes A)^L$. The claims of the theorem follow.
\end{proof}

\section{Decomposition theorem for the jets of vector fields}

In this section, we apply the Taylor formula to establish a decomposition theorem for the jets of vector fields on $G$. Since $V = A \otimes \g^L$, we get that 
$$A \widehat\# V = (A \widehat{\mathop\otimes\limits_\Delta} A ) \otimes \g^L.$$
Next, applying (\ref{transform}), we obtain
$$A \widehat\# V \cong 
(A \mathop{\widehat\otimes}\limits_{1 \otimes \m_e} A) \otimes \g^L.$$
Taking into account the decomposition $A = K \cdot 1 \oplus \m_e$, we get 
$$A \widehat\# V \cong (A \otimes 1) \otimes \g^L \oplus  
(A \mathop{\widehat\otimes}\limits_{1 \otimes \m_e} \m_e) \otimes \g^L.$$
We claim that the first summand on the right-hand side forms a Lie subalgebra, while the second summand is an ideal. Let us compute the Lie brackets.

\begin{lemma}
\label{V}
The subspace $A \# \g^L$ in $A \# V$ forms a subalgebra, which is isomorphic to $V$.
\end{lemma}
\begin{proof}
By comparing formula (\ref{smash}) for the Lie brackets in the smash product with the formula for the Lie bracket in $V$, we see that the map $A \# \g^L \rightarrow V$, $f \# x \mapsto fx$, for $f\in A$, $x \in \g^L$, is an isomorphism of Lie algebras.
\end{proof}

We point out that the subalgebra $A \# \g^L$ is distinct from the natural subalgebra of vector fields $1 \# V$ of the Lie algebra of algebraic jets of vector fields $A \# V$.

Denote by $\LL_+$ the subalgebra of vector fields vanishing at $e$: 
$\LL_+ = \m_e \otimes \g^L$. Since $\m_e^N \otimes \g^L$ is an ideal in $\LL_+$, we can take the completion: 
$$\hL_+ = \varprojlim_N (\m_e \otimes \g^L) / (\m_e^N \otimes \g^L).$$

We will also need the completion 
$$A \widehat\otimes \LL_+ = \varprojlim_N A \otimes \left( (\m_e \otimes \g^L) / (\m_e^N \otimes \g^L) \right).$$
We point out that tensor products do not commute with inverse limits in general.

Let $L = \Der K[t_1, \ldots, t_n]$ and 
$\widehat{L} = \Der K[[t_1, \ldots, t_n]]$ with $n = \dim G$. Let $L_+ = \m_0 L$ and
$\widehat{L}_+ = \m \widehat{L}$, 
where $\m_0$ is the ideal of functions vanishing at $0$ and
$\m$ is the unique maximal ideal of $K[[t_1, \ldots, t_n]]$.
By \cite{BFN} (Lemma 6), $\hL_+ \cong \widehat{L}_+$.  
Also, for every $N$, 
$$\hL_+ / \m_e^N \hL_+ \,\cong\, \widehat{L}_+ / \m^N \widehat{L}_+ \,\cong\, L_+ / \m_0^N L_+.$$

\begin{lemma}
\label{AL}
The subspace $(A \mathop{\widehat\otimes}\limits_{1 \otimes \m_e} \psi (\m_e)) \otimes \g^L$ in $A \widehat\# V$ is a Lie subalgebra isomorphic to $A \widehat\otimes \LL_+$
with Lie bracket 
$$[f_1 \otimes \tau_1, f_2 \otimes \tau_2] = f_1 f_2 \otimes [\tau_1, \tau_2],
\ \ \text{for } f_1, f_2 \in A, \tau_1, \tau_2 \in \LL_+,$$ 
extended to the completion by continuity.
\end{lemma}
\begin{proof}
We need to compute the Lie bracket $[(f_1 \otimes 1) \psi(h_1) \otimes x^L,  (f_2 \otimes 1) \psi(h_2) \otimes y^L]$ in $A \# V$. Since $\psi(h_1), \psi(h_2) \in \Delta$, the terms with $\mu$ in (\ref{doublesmash}) vanish, and we get
\begin{align*}
(f_1 f_2 \otimes 1) \psi(h_1) (1 \otimes x^L) \psi(h_2) \otimes y^L
&- (f_1 f_2 \otimes 1) \psi(h_2) (1 \otimes y^L) \psi(h_1) \otimes x^L \\
&+ (f_1 f_2 \otimes 1) \psi(h_1) \psi(h_2) \otimes [x,y]^L \\
\end{align*}
Applying Corollary \ref{transfer}, we get
\begin{equation*}
(f_1 f_2 \otimes 1) \psi(h_1 x^L(h_2)) \otimes y^L
- (f_1 f_2 \otimes 1) \psi(h_2 y^L (h_1)) \otimes x^L 
+ (f_1 f_2 \otimes 1) \psi(h_1 h_2) \otimes [x,y]^L .
\end{equation*}
This matches the Lie bracket in $A \otimes \LL_+$.
\end{proof}

Note that $\m_e$ is a $G$-module under the conjugation action $C^*$. This follows from the fact that the identity element is fixed by conjugation. This $G$-action induces a $\g$-module structure on $\m_e$, where $x \in \g$ acts as vector field $x^C = x^L + x^R$. The fact that $x^L$ and $x^R$ evaluate to opposite vectors at $e$, gives us an alternative proof of invariance of $\m_e$ under this action of $\g$.   

Now we are ready to prove the main result of the paper:

\begin{theorem} [Decomposition Theorem]
\label{decomposition}
Let $G$ be a connected affine algebraic group. Then $A \widehat\# V$ decomposes into a semidirect product:
$$A \widehat\# V \cong V \ltimes A \widehat\otimes \LL_+$$
with the following action of $V$ on the ideal $A \widehat\otimes \LL_+$:
$$[f_1 \otimes x^L, f_2 \otimes h \otimes y^L ] = f_1 x^L(f_2)   \otimes h \otimes y^L 
+ f_1 f_2 \otimes x^C(h) \otimes y^L + f_1 f_2 \otimes h \otimes [x,y]^L,$$
extended to the completion by continuity.
In this isomorphism, the embedding maps $V \rightarrow A \# V$ and $A \otimes \LL_+ \rightarrow  A \# V$ are given by
\begin{align*}
f \otimes x^L &\mapsto f \# x^L, \\
f \otimes h \otimes y^L &\mapsto (f \# 1) \psi(h)  y^L,
\end{align*}
extended to completions by continuity.
\end{theorem}
\begin{remark}
The semidirect product $V \ltimes A \widehat\otimes \LL_+$ has a subalgebra $\g \ltimes \hL_+$ with the Lie bracket
$$[x, h \otimes y] = x^C (h) \otimes y + h \otimes [x, y], \ \ \text{for \,} h \in \m_e, x,y\in\g.$$
The Lie algebra $V \ltimes A \widehat\otimes \LL_+$ is a completed smash product of $\g \ltimes \hL_+$ with $A$, where $x \in \g$ acts on $A$ via $x^L$, and $\hL_+$ acts on $A$ trivially.
\end{remark}
\begin{proof}
We have already seen in Lemmas \ref{V} and \ref{AL}
that the maps $V \rightarrow A \# V$ and $A \otimes \LL_+ \rightarrow  A \# V$ are embeddings
of Lie algebras.
We need to compute the Lie bracket of the images of $f_1 \otimes x^L$ and $f_2 \otimes h \otimes y^L$ in $A \# V$. For this, we use (\ref{doublesmash}):
\begin{align*}
&[(f_1 \# 1) x^L, (f_2 \# 1) \psi(h) y^L]  \\
&= (f_1 \# 1) (x^L (f_2) \# 1) \psi(h) y^L 
+  (f_1 f_2 \# 1)(x^L \otimes 1) \psi(h) y^L 
+  (f_1 f_2 \# 1)(1 \otimes x^L) \psi(h) y^L\\
&- (f_2 \mu(\psi(h)) \# 1) (y^L (f_1) \# 1) x^L
+ (f_1 f_2 \# 1)  \psi(h) [x,y]^L.
\end{align*}
Since $\psi(h) \in \Delta$ for $h\in\m_e$, we get that $\mu(\psi(h)) = 0$. By Corollary \ref{transfer}, 
$(x^L \otimes 1 + 1\otimes x^L) \psi(h) = \psi(x^C (h))$. The above expression simplifies to
\begin{align*}
(f_1 x^L (f_2) \# 1) \psi(h) y^L 
+  (f_1 f_2 \# 1) \psi(x^C(h)) y^L 
+ (f_1 f_2 \# 1)  \psi(h) [x,y]^L.
\end{align*}
This matches the Lie bracket in the semidirect product $V \ltimes A \widehat\otimes \LL_+$.

We point out that the embedding $V \ltimes A \otimes \LL_+ \rightarrow A \# V$ respects the filtrations: the image of $\m_e^k \LL_+$ is in $\Delta^{k+1} \otimes V$. This means that this embedding extends by continuity to a homomorphism of completed Lie algebras, which becomes an isomorphism by the Taylor formula.
\end{proof}

The most important class of $AV$-modules is the category of \emph{differentiable} $AV$-modules. Let us recall the definition.

\begin{definition}
An $AV$-module $M$ is called $N$-differentiable if it is annihilated by operators
$$\sum_{k=0}^N (-1)^k {N \choose k} f^{N-k} \# f^k \tau$$
for all $f \in A$, $\tau \in V$. 
\end{definition}
Equivalently, an $AV$-module is $N$-differentiable if it is annihilated by $\Delta^N \otimes V$.
$D$-modules are precisely 1-differentiable $AV$-modules.
\begin{corollary}
\label{Ndiff}
$N$-differentiable $AV$-modules on an affine algebraic group $G$ are precisely modules for the associative algebra $A \# U(\g \ltimes (\LL_+ / \m_e^{N-1} \LL_+))$, which is the smash product of a Hopf algebra with its module $A$. Here $\g$ acts on $A$ by left-invariant derivations, whereas the action of $\LL_+$ on $A$ is trivial. 
\end{corollary}
\begin{proof}
$N$-differentiable $AV$-modules are strong modules for the Rinehart pair \break
$(A, \, A\#V / (\Delta^N \otimes V))$. Since we are considering here a truncation of a filtration, we have an isomorphism $A\#V / (\Delta^N \otimes V) \cong A\widehat\#V / (\widehat\Delta^N \otimes V)$.
Applying the Decomposition Theorem, we get 
$A\widehat\#V / (\widehat\Delta^N \otimes V) \cong V \ltimes (A \widehat\otimes \LL_+)/(A \widehat\otimes \m_e^{N-1} \LL_+)$. Again, since we truncate the filtration, we have 
$(A \widehat\otimes \LL_+)/(A \widehat\otimes \m_e^{N-1} \LL_+) \cong A \otimes (\LL_+ / \m_e^{N-1}\LL_+)$. We conclude that $A\#V / (\Delta^N \otimes V)$ is isomorphic to the smash product
$A \# (\g \ltimes (\LL_+ / \m_e^{N-1}\LL_+))$. However,
\begin{align*}
U_\text{strong} (A, \, A \# (\g \ltimes (\LL_+ / \m_e^{N-1}\LL_+))) 
&\cong U_\text{weak} (A, \, \g \ltimes (\LL_+ / \m_e^{N-1}\LL_+)) \\
&\cong A \# U( \g \ltimes (\LL_+ / \m_e^{N-1}\LL_+)),
\end{align*}
which establishes the claim.
\end{proof}

\section{Free $AV$-modules of finite rank and the Lax equation}

In this final section, we will discuss a few examples of $AV$-modules on affine algebraic groups. We say that an $AV$-module has finite rank (resp. free) if it is finitely generated (resp. free) as an $A$-module. 

The first natural class of $AV$-modules has a geometric origin. These are tensor modules. The space of vector fields $V$ on $G$ is a free $AV$-module of finite rank, $V \cong A \otimes \g$. Its dual is the space of 1-forms,
$\Omega^1 = \Hom_A (V, A)$, which is also a free $A$-module, $\Omega^1 \cong A \otimes \g^*$. 
By taking the tensor product over $A$ of $k$ copies of $V$ and $m$ copies of $\Omega^1$, we obtain the $AV$-module of $(k,m)$-tensor fields, which will again be a free module of finite rank. All tensor modules are 2-differentiable.

Let us generalize the construction of $(k, m)$-tensor modules. By Corollary \ref{Ndiff}, every 2-differentiable $AV$-module is a module for $\g \ltimes (\LL_+/\m_e \LL_+)$. Note that $\LL_+/\m_e \LL_+ \cong L_+ / \m_0 L_+ \cong \gl(\g)$.

\begin{proposition}
In the semidirect product $\g \ltimes (\LL_+/\m_e \LL_+) \cong \g \ltimes \gl(\g)$, the action of $\g$ on 
$\gl(\g) \cong \g \otimes \g^*$ is the tensor product of the adjoint and coadjoint representations. 
\end{proposition}
\begin{proof}
In the semidirect product, $\g^L$ acts on $(\m_e /\m_e^2) \otimes \g^L$ via the tensor product of the conjugation action and the adjoint action. Consider the pairing
$$\m_e / \m_e^2 \times \g^L \rightarrow A /\m_e \cong K,$$
given by $\left<f, x^L\right> = x^L(f)(e)$. Let us check that this is an invariant pairing of $\g$-modules, where $x \in \g$ acts on  $\m_e / \m_e^2$ via $x^C$ and on $\g^L$ via $\ad(x^L)$. Indeed for $y \in \g$,
\begin{align*}
\left<y^C (f), x^L\right> &= x^L (y^C (f))(e) = x^L (y^L (f))(e) + x^L (y^R(f))(e) \\
\left<f, \ad(y^L) x^L\right> &= [y^L, x^L] (f)(e) =  - x^L (y^L (f))(e) +  y^L (x^L (f))(e).
\end{align*}
Since $x^L$ and $y^R$ commute, we have 
$$\left<y^C (f), x^L\right> + \left<f, \ad(y^L) x^L\right> = y^C (x^L(f))(e).$$
However, the vector field $y^C$ vanishes at $e$, which means that $\left<y^C (f), x^L\right> + \left<f, \ad(y^L) x^L\right> = 0$. The space $\m_e /\m_e^2$ is identified with $\g^*$ via the above pairing, and the invariance property means that it is isomorphic to the coadjoint $\g$-module. Hence, $\LL_+/\m_e \LL_+$ is isomorphic to 
$\g^* \otimes \g$ as a $\g$-module.
\end{proof}

\begin{proposition}
\label{splits}
Consider the action of $\g$ on $\gl(\g)$ via the tensor product of the adjoint and coadjoint actions. The corresponding semidirect product splits, $\g \ltimes \gl(\g) \cong \g \oplus \gl(\g)$ as Lie algebras.
\end{proposition}
\begin{proof}
The action of $\g$ on $\gl(\g)$ gives a homomorphism $\rho: \ \g \rightarrow \Der (\gl(\g))$. It is easy to see that the derived subalgebra and the center of a Lie algebra are invariant under every derivation of that Lie algebra. This means that the subspaces $K I$ and $\sll(\g)$ are invariant under $\rho(x)$ for every $x \in \g$.
We claim that $\rho(x) I = 0$. Since $K I$ is an invariant subspace, it is sufficient to show that $\tr (\rho(x) I) = 0$.
Choose dual bases $\{y_i\}$ of $\g$ and $\{y_i^*\}$ of $\g^*$. Then $I = \sum_i y_i \otimes y_i^*$, 
$\rho(x) I = \sum_i \ad(x) y_i \otimes y_i^* + y_i \otimes \ad^*(x) y_i^*$ and
$$\tr(\rho(x) I ) = \sum_i y_i^* ([x,y_i]) + \ad^*(x) y_i^*(y_i)
= \sum_i y_i^* ([x,y_i]) - y_i^* ([x,y_i]) = 0.$$
We conclude that $\rho: \ \g \rightarrow \Der (\sll(\g))$. However, every derivation of $\sll(\g)$ is inner,
$\Der (\sll(\g)) \cong \sll(\g)$. Composing these maps, we get a homomorphism $\widetilde\rho: \g \rightarrow \sll(\g)$ that realizes the action by inner derivations. It is easy to see that elements $x - \widetilde\rho(x)$ span a subalgebra of $\g \ltimes \gl(\g)$ which is isomorphic to $\g$ and commutes with $\gl(\g)$. This proves that the semidirect product splits.
\end{proof}

\begin{corollary}
A 2-differentiable $AV$-module on an affine algebraic group $G$ is precisely a $D$-module on $G$ equipped with an additional action of $\gl_n$, $n = \dim G$, where the action of $\gl_n$ commutes with the action of $D$.
\end{corollary}

Let us present explicit formulas for the $AV$-module action in the setting of the above corollary.

Given $f \in A$, $x \in \g$, we view
$f \ad(x)$ as an element of $A \otimes \gl(\g)$. In the Lie algebra $A \otimes \gl(\g)$, $\gl(\g)$ acts on $A$ trivially. 

For $f \in A$, the differential $df$ belongs to $A \otimes \g^*$, with $df (y) = y^L(f)$. 
Then $x df \in A \otimes \g \otimes \g^* = A \otimes \gl(\g)$.


\begin{theorem}
[ChatGPT]
\label{twodiff}
Let $G$ be an affine algebraic group. Let $M$ be a $D$-module on $G$, and simultaneously a $\gl(\g)$-module, with commuting actions. Then $M$ admits a structure of an $AV$-module with $A$-action defined by the $D$-module structure and the following action of $V$:
$$\rho(f x^L) = f x^L + f \ad(x) - x df, \ \ \text{for } f \in A, x \in \g.$$
The right-hand side is viewed as an operator on $M$, where the first term comes from the $D$-module action, and the other two terms correspond to the action of $A \otimes \gl(\g)$.
\end{theorem}
\begin{proof}
We point out that the action of $\gl(\g)$ on $M$ is $A$-linear, since the action of $\gl(\g)$ commutes with the action of $D$, and hence the action of $A$. We need to show that the formula defines a representation of $V$ on $M$. For this, we record the Lie brackets involving the operators appearing on the right-hand side.
\begin{align*}
[f x^L, h y^L] &= f x^L(h) y^L - h y^L(f) x^L + fh [x,y]^L, \\
[f x^L, h \,\ad(y)] &= f x^L(h) \,\ad(y), \\
[f \,\ad(x), h \,\ad(y)] &= fh \, \ad([x,y]), \\
[x \,df, y \,dh] &= y^L(f) x \,dh - x^L(h) y \,df, \\
[f x^L + f \,\ad(x), y \,dh] &= f [x,y] \,dh + f y \,d(x^L(h)).
\end{align*}
These relations are straightforward to verify, and the fact that $\rho$ is a representation of $V$ easily follows.
\end{proof}

In particular, if $M$ is a $D$-module on $G$ and $W$ is a $\gl_n$-module, then $M \otimes W$ admits the structure of a 2-differentiable $AV$-module on $G$. 

Let us now discuss free $AV$-modules of finite rank. By the theorem of Bouaziz-Rocha \cite{BR}, every finitely generated $AV$-module is differentiable, and hence by Corollary \ref{Ndiff}, admits the action of 
$A \# U(\g \ltimes (\LL_+ / \m_e^{N-1} \LL_+))$ for some $N$. This algebra contains a subalgebra 
$A \# U(\g)$, which is the algebra $D$ of differential operators on $G$. Hence $N$-differentiable $AV$-modules are $D$-modules with an additional action of the Lie algebra $\LL_+ / \m_e^{N-1} \LL_+$, compatible with the $D$-module structure.

Let $M$ be a free $D$-module of finite rank, $M = A \otimes W$, $\dim W < \infty$. The $D$-module structure is determined by the action of $\g^L$ on $W$: $\g^L \times W \rightarrow A \otimes W$.
Such a map is encoded by a 1-form $\omega$ on $G$ with values in $\End(W)$. Then the action of $V$ on $M$ is given by
$$(f x^L) (h w) = f x^L(h)  w + fh \,\omega(x^L)(w), \ \ \text{for \,} f,h\in A, \,x\in \g, \,w\in W.$$ 
In order for this to yield a Lie algebra action, the form $\omega$ must satisfy the condition
\begin{equation}
\label{domega}
d \omega + \omega \wedge \omega = 0,
\end{equation}
or equivalently,
$$\omega([x,y]^L) = x^L (\omega(y^L)) -  y^L (\omega(x^L)) + [\omega(x^L), \omega(y^L)].$$ 

Let $\dim W = m$. A change of basis in $A \otimes W$ is described by the group $GL_m(A)$. Under a change of basis given by $g \in GL_m(A)$, the 1-form $\omega$ transforms as \cite{Bo}: 
$$\omega^\prime = g^{-1} \omega g + g^{-1} dg.$$

Let us consider the simplest case, $m = 1$. The condition on $\omega$ in this case is just $d \omega = 0$, i.e., the form must be closed. In particular, we may take an exact form, $\omega = dH$, $H \in A$, and the action of $V$ on this $D$-module will be given by
$$(f x^L) (h w) = f x^L(h) w + fh x^L(H) w.$$
Note that $GL_1 (A) = A^*$. If the units of $A$ are scalars, two such $D$-modules are non-isomorphic unless $H_1 - H_2$ is a constant function \cite{BF}.

Now let us see how a free rank-one $ D$-module can be turned into an $AV$-module. For this, we need to define a compatible action of $\LL_+ / \m_e^{N-1} \LL_+$. Since $\End(W)$ is abelian, the derived subalgebra
$\left( \LL_+ / \m_e^{N-1} \LL_+ \right)^{(1)}$ must act trivially. The quotient by the derived subalgebra is 1-dimensional and is spanned by the identity operator $I$ in $\gl(\g) = \LL_+ / \m_e \LL_+$. Assume that $I$ acts as multiplication by $\lambda \in A$. Since the action of $V$ must commute with the action of $I$, we see that $x^L(\lambda) = 0$ for all $x \in \g$, so $\lambda$ must be a scalar in $K$.
The action of $V$ on the corresponding $AV$-module is 
$$(f x^L) (h w) = f x^L(h) w + fh\, \omega(x^L) w + \lambda h x^L(f) w.$$

Let us generalize the above construction of an $AV$-module to higher rank. Let $B \in M_m(K[t_1,\ldots, t_k])$ satisfy the condition
\begin{equation}
\label{commute}
\left[ \frac{\partial B}{\partial t_i},  \frac{\partial B}{\partial t_j} \right] = 0 \ \ \text{for all \,} i,j.
\end{equation}
Fix $H_1, \ldots, H_k \in A$. Then $\omega = dB(H_1, \ldots, H_k)$ satisfies (\ref{domega}) and yields a $D$-module structure on $A \otimes K^m$:
\begin{equation}
\label{Dmod}
(f x^L)(hw) = f x^L(h) w + \sum_{j=1}^k fh x^L(H_j) \frac{\partial B}{\partial t_j} (H_1, \ldots, H_k) w.
\end{equation}

\begin{theorem}
Let $G$ be an affine algebraic group.
Let $B, L \in M_m(K[t_1,\ldots, t_k])$ satisfy (\ref{commute}) and
\begin{equation}
\label{Lax}
\frac{\partial L}{\partial t_i} = \left[ L,   \frac{\partial B}{\partial t_i} \right]
\ \ \text{for \,} i=1,\ldots,k.
\end{equation}
Let $H_1, \ldots, H_k \in A$. Then we have an $AV$-module structure on $A \otimes K^m$, where the action of $V$ is given by
\begin{equation}
\label{AVaction}
(f x^L)(hw) = f x^L(h) w + \sum_{j=1}^k fh x^L(H_j) \frac{\partial B}{\partial t_j} (H_1, \ldots, H_k) w + x^L(f) h L(H_1, \ldots, H_k) w.
\end{equation}
\end{theorem}
\begin{proof}
To construct a differentiable $AV$-module on $G$, we need to take a $D$-module on $G$ and define on it a compatible $A$-linear action of $\LL_+ / \m_e^{N-1} \LL_+$.
For the construction presented in this theorem, we take a $D$-module given by (\ref{Dmod}), and we define the action of $\LL_+ / \m_e^{N-1} \LL_+$, where the derived subalgebra is acting trivially and the identity element $I \in \gl(\g)$ acts by 
$L(H_1, \ldots, H_k)$. Let us verify that these actions are compatible. The compatibility condition for this setting is that the action of $I$ must commute with the $D$-module action (see the proof of Proposition \ref{splits}). Indeed, using (\ref{Lax}), we get
\begin{align*}
&\left[ x^L + \sum_{j=1}^k x^L(H_j) \frac{\partial B}{\partial t_j} (H_1, \ldots, H_k), 
\ L(H_1, \ldots, H_k) \right] = \\
&\sum_{j=1}^k x^L(H_j) \frac{\partial L}{\partial t_j} (H_1, \ldots, H_k)
+ \sum_{j=1}^k x^L(H_j) \left[ \frac{\partial B}{\partial t_j} (H_1, \ldots, H_k), 
\ L(H_1, \ldots, H_k) \right] = 0.
\end{align*}
Applying Theorem \ref{decomposition} to combine the two actions, we recover 
(\ref{AVaction}).
\end{proof}

We point out that when $k=1$, condition (\ref{commute}) holds trivially, while (\ref{Lax}) becomes
$$\frac{dL}{dt} = \left[ L,   \frac{dB}{dt} \right].$$
This is known as the Lax equation in the theory of integrable systems. Many integrable systems may be expressed through the Lax equation. Finite-dimensional integrable systems of classical mechanics are written in Lax form in the context of algebraic groups. The Kovalevskaya top provides one such example \cite{Ts}.
System (\ref{Lax}) is a generalized Lax equation with commuting flows.
Thus, the construction reveals a connection between the Lax equation and the theory of $AV$-modules.

\begin{conjecture}
Let $A$ be the algebra of regular functions on a connected simply-connected affine algebraic group. Every finitely generated projective $A$-module is free.
\end{conjecture}

\end{document}